\documentclass[12pt]{article}
\usepackage{amssymb}
\usepackage{latexsym}
\usepackage{amscd}

\newtheorem{definition}{Definition}[section]
\newtheorem{theorem}[definition]{Theorem}
\newtheorem{lemma}[definition]{Lemma}

\newtheorem{remark}[definition]{Remark}
\newtheorem{example}[definition]{Example}

\newtheorem{note}[definition]{Note}
\newtheorem{assumption}[definition]{Assumption}
\newtheorem{proposition}[definition]{Proposition}

\def\R{\mathbb R}
\def\C{\mathbb C}

\def\Z{\mathbb Z}

\begin{document}

\title{\bf Distance-regular 
graphs of 
$q$-Racah\\ type and the $q$-tetrahedron algebra}

\author{
Tatsuro Ito{\footnote{
Department of Computational Science,
Faculty of Science,
Kanazawa University,
Kakuma-machi,
Kanazawa 920-1192, Japan
}}
{\footnote{Supported in part by JSPS grant 18340022}} 
$\;$ and
Paul Terwilliger{\footnote{
Department of Mathematics, University of
Wisconsin, 480 Lincoln Drive, Madison WI 53706-1388 USA}
}}
\date{}
%to get date printout, comment out above line

\maketitle

\centerline{\large In Memory of Donald Higman} 

\begin{abstract}
In this paper we discuss a relationship between the
following two algebras:
(i)
the subconstituent algebra $T$ of  a 
distance-regular graph that has $q$-Racah type;
(ii) the $q$-tetrahedron algebra 
$\boxtimes_q$ which is a 
$q$-deformation of 
the three-point 
${\mathfrak{sl}}_2$ loop algebra.
Assuming that every irreducible $T$-module is thin, we 
display an algebra homomorphism from $\boxtimes_q$ into $T$ and
show that $T$ is generated by the image
together with the center $Z(T)$.

\medskip
\noindent
{\bf Keywords}. Tetrahedron algebra,
quantum affine algebra,
distance-regular graph,
$Q$-polynomial.
 \hfil\break
\noindent {\bf 2000 Mathematics Subject Classification}. 
Primary: 05E30. Secondary: 05E35; 17B37.
 \end{abstract}

\section{Introduction}
In \cite{HT}
B. Hartwig and the second author
gave a presentation of the three-point $\mathfrak{sl}_2$ loop
algebra via generators and relations. To obtain this
presentation they defined a Lie algebra $\boxtimes$ by
generators and relations, and displayed an isomorphism
 from $\boxtimes$ to
the three-point $\mathfrak{sl}_2$ loop algebra. The algebra
$\boxtimes$
is called the tetrahedron algebra \cite[Definition 1.1]{HT}.
In \cite{qtet} we introduced a
$q$-deformation $\boxtimes_q$ of 
$\boxtimes$ called the
$q$-tetrahedron algebra.
In \cite{qtet} and \cite{qinv}
we described the finite-dimensional irreducible
$\boxtimes_q$-modules.
In \cite[Section 4]{qtetanddrg} we displayed four homomorphisms
into $\boxtimes_q$ from the quantum affine algebra
$U_q({\widehat{\mathfrak{sl}}_2})$.
In
\cite[Section 12]{qtetanddrg} we found a homomorphism
from
$\boxtimes_q$ into the subconstituent algebra of a distance-regular
graph that is self-dual with classical parameters.
In the present paper we do something similar 
for a 
 distance-regular graph said to have $q$-Racah type.
This type is described as follows.
Let $\Gamma$ denote a distance-regular graph with diameter
$D\geq 3$ (See Section 4 for formal definitions).
We say that $\Gamma$ has {\it $q$-Racah type} whenever
$\Gamma$ has a $Q$-polynomial structure with eigenvalue
sequence
$\lbrace \theta_i \rbrace_{i=0}^D$
and dual eigenvalue sequence
$\lbrace \theta^*_i \rbrace_{i=0}^D$ that satisfy
\begin{eqnarray*}
\theta_i &=& \eta + uq^{2i-D}+vq^{D-2i} \qquad (0 \leq i \leq D),
\\
\theta^*_i &=& \eta^* + u^*q^{2i-D}+v^*q^{D-2i} \qquad (0 \leq i \leq D),
\end{eqnarray*}
 where $q,u,v,u^*,v^*$ are nonzero
and $q^{2i}\not=1$ for
$1 \leq i \leq D$.
Assume $\Gamma$ has $q$-Racah type.
Fix a vertex $x$ of $\Gamma$ and let $T=T(x)$ denote the corresponding
subconstituent algebra
\cite[Definition 3.3]{terwSub1}.
 Recall that $T$ is generated by the
adjacency matrix $A$ and the dual adjacency matrix $A^*=A^*(x)$
\cite[Definition 3.10]{terwSub1}.
An irreducible $T$-module $W$ is called {\it thin}
whenever the intersection of $W$ with each eigenspace of
$A$ and each eigenspace of $A^*$ has dimension at most 1
\cite[Definition 3.5]{terwSub1}.  
Assuming
each irreducible $T$-module is
thin, we  
display 
invertible central elements $\Phi, \Psi$ of $T$ and
a homomorphism $\vartheta : \boxtimes_q \to T$
such that
\begin{eqnarray*}
A&=&
\eta I + u \Phi \Psi^{-1}\vartheta(x_{01})+v\Psi \Phi^{-1}\vartheta(x_{12}),
\\
A^*
&=&\eta^* I + u^* \Phi \Psi \vartheta(x_{23})+v\Psi^{-1} \Phi^{-1}
\vartheta(x_{30}),
\end{eqnarray*}
 where the $x_{ij}$ are the standard generators of $\boxtimes_q$.
It follows that $T$ is generated by the image
$\vartheta(\boxtimes_q)$ 
together with $\Phi, \Psi$. In particular $T$ is generated by
$\vartheta(\boxtimes_q)$ together with the center $Z(T)$.

\medskip
\noindent 
This paper is organized as follows. In Section 2 we
recall the definition of $\boxtimes_q$.
In Section 3 we describe how $\boxtimes_q$ is related to
$U_q({\widehat{\mathfrak{sl}}}_2)$. 
In Section 4 we recall the basic theory of a distance-regular graph
$\Gamma $, focussing on the $Q$-polynomial
property and the subconstituent algebra.
In Section 5 we discuss the split decomposition of $\Gamma$.
In Section 6 we give our main results.

\medskip
\noindent Throughout the paper $\C$ denotes the field of
complex numbers.

\section{The $q$-tetrahedron algebra $\boxtimes_q$}

\noindent In this section we recall the 
$q$-tetrahedron algebra. 
We fix a nonzero scalar $q \in \C$ such that $q^2\not=1$ and define
\begin{eqnarray*}
\lbrack n \rbrack_q = \frac{q^n-q^{-n}}{q-q^{-1}},
\qquad \qquad n = 0,1,2,\ldots 
\label{eq:nbrack}
\end{eqnarray*}
We let $\Z_4 = \Z/4\Z$ denote the cyclic group of order 4.

\begin{definition} \rm \cite[Definition 10.1]{qtet}
\label{def:qtet}
Let $\boxtimes_q$ denote the unital associative $\C$-algebra that has
generators 
\begin{eqnarray*}
\lbrace x_{ij}\;|\; i,j \in \Z_4,\;j-i=1 \;\mbox{or} \;j-i=2\rbrace
\end{eqnarray*}
and the following relations:
\begin{enumerate}
\item For $i,j\in \Z_4$ such that $j-i=2$,
\begin{eqnarray*}
x_{ij}x_{ji} = 1.
\label{eq:qrel0}
\end{eqnarray*}
\item For $h,i,j\in \Z_4$ such that the pair $(i-h,j-i)$ is one of
$(1,1), (1,2), (2,1)$,
\begin{eqnarray*}
\frac{qx_{hi}x_{ij}-q^{-1}x_{ij}x_{hi}}{q-q^{-1}}=1.
\label{eq:qrel1}
\end{eqnarray*}
\item For $h,i,j,k\in \Z_4$ such that $i-h=j-i=k-j=1$,
\begin{eqnarray}
\label{eq:qserre}
x_{hi}^3x_{jk} -
\lbrack 3 \rbrack_q
x_{hi}^2x_{jk}x_{hi} +
\lbrack 3 \rbrack_q
x_{hi}x_{jk}x_{hi}^2- 
x_{jk}x_{hi}^3=0. 
\end{eqnarray}
\end{enumerate}
We call $\boxtimes_q$ the
{\it $q$-tetrahedron algebra} or ``$q$-tet'' for short.
We refer to the $x_{ij}$ as the {\it standard generators}
for $\boxtimes_q$.
\end{definition}
\begin{note}\rm
The equations (\ref{eq:qserre}) are the cubic $q$-Serre relations
\cite[p.~10]{lusztig}. 
\end{note}
\noindent
We make some observations.

\begin{lemma}
\label{lem:rho}
{\rm \cite[Lemma 6.3]{qtet} }
There exists a $\C$-algebra automorphism $\varrho$ of $\boxtimes_q$
that sends each generator $x_{ij}$ to $x_{i+1,j+1}$.
Moreover 
 $\varrho^4=1$.
\end{lemma}

\begin{lemma}
\label{lem:flip}
{\rm \cite[Lemma 6.5]{qtet}}
There exists a $\C$-algebra automorphism of 
 $\boxtimes_q$
that sends each generator $x_{ij}$ to $-x_{ij}$.
%This automorphism has order 2.
\end{lemma}

\section{The quantum affine algebra
$U_q(\widehat{ \mathfrak{sl}}_2)$
}

\noindent In this section we consider how
$\boxtimes_q$ is related to the
quantum affine algebra
$U_q(\widehat{ \mathfrak{sl}}_2)$.
We start with a definition.

\begin{definition} 
\label{def:qauq}
\rm
\cite[p.~266]{charp} 
The quantum affine algebra
$U_q(\widehat{ \mathfrak{sl}}_2)$
is the unital associative $\C$-algebra 
with
generators $K^{\pm 1}_i$,
 $e^{\pm}_i$,
$i\in \lbrace 0,1\rbrace $
and the following relations:
\begin{eqnarray*}
K_iK^{-1}_i &=& K^{-1}_iK_i=1,
\label{eq:qauq1}
\\
K_0K_1&=& K_1K_0,
\label{eq:qauq2}
\\
K_ie^{\pm}_iK^{-1}_i &=& q^{{\pm}2}e^{\pm}_i,
\label{eq:qauq3}
\\
K_ie^{\pm}_jK^{-1}_i &=& q^{{\mp}2}e^{\pm}_j, \qquad i\not=j,
\label{eq:qauq4}
\\
\lbrack e^+_i, e^-_i\rbrack &=& {{K_i-K^{-1}_i}\over {q-q^{-1}}},
\label{eq:qauq5}
\\
\lbrack e^{\pm}_0, e^{\mp}_1\rbrack &=& 0,
\label{eq:qauq6}
\end{eqnarray*}
\begin{eqnarray*}
(e^{\pm}_i)^3e^{\pm}_j -  
\lbrack 3 \rbrack_q (e^{\pm}_i)^2e^{\pm}_j e^{\pm}_i 
+\lbrack 3 \rbrack_q e^{\pm}_ie^{\pm}_j (e^{\pm}_i)^2 - 
e^{\pm}_j (e^{\pm}_i)^3 =0, \qquad i\not=j.
\label{eq:qauq7}
\end{eqnarray*}
%We call $e^{\pm}_i$, $K_i^{{\pm}1}$, $i\in \lbrace 0,1\rbrace $
%the {\it Chevalley generators} for
%$U_q({\widehat{sl}}_2)$.
\end{definition}

\noindent The following presentation of
$U_q(\widehat{ \mathfrak{sl}}_2)$
 will be useful.

\begin{proposition}
\label{thm:qa2} 
{\rm (\cite[Theorem 2.1]{tdanduq},
\cite{equit2})}
The quantum affine algebra
$U_q(\widehat{ \mathfrak{sl}}_2)$
 is isomorphic to
the unital associative $\C$-algebra 
with
generators $x_i^{\pm 1}$, $y_i$, $z_i$, $i\in \lbrace 0,1\rbrace $
and the following relations:
\begin{eqnarray*}
x_ix^{-1}_i = x^{-1}_ix_i &=&1,\\
x_0x_1 \;\;\mbox{is central},
\label{eq:qabuq2}
\\
\frac{q x_iy_i-q^{-1}y_ix_i}{q-q^{-1}} &=& 1,
\label{eq:qabuq3}
\\
\frac{q y_iz_i-q^{-1}z_iy_i}{q-q^{-1}} &=& 1,
\label{eq:qabuq4}
\\
\frac{q z_ix_i-q^{-1}x_iz_i}{q-q^{-1}} &=& 1,
\label{eq:qabuq5}
\\
\frac{q z_iy_j-q^{-1}y_jz_i}{q-q^{-1}} &=&x^{-1}_0x^{-1}_1,
\qquad i\not=j,
\label{eq:qabuq6}
\end{eqnarray*}
\begin{eqnarray*}
y_i^3y_j -  
\lbrack 3 \rbrack_q y_i^2y_j y_i 
+\lbrack 3 \rbrack_q y_iy_j y_i^2 - 
y_j y_i^3 =0, \qquad i\not=j,
\label{eq:qabuq7}
\\
z_i^3z_j -  
\lbrack 3 \rbrack_q z_i^2z_j z_i 
+\lbrack 3 \rbrack_q z_iz_j z_i^2 - 
z_j z_i^3 =0, \qquad i\not=j.
\label{eq:qabuq8*}
\end{eqnarray*}
An isomorphism with the presentation in Definition
\ref{def:qauq} is given by:
\begin{eqnarray*}
\label{eq:qaiso1}
x^{\pm 1}_i &\mapsto & K^{\pm 1}_i,\\
\label{eq:qaiso2}
y_i &\mapsto & K^{-1}_i+e^{-}_i, \\
\label{eq:qaiso3}
z_i &\mapsto & K^{-1}_i-K^{-1}_ie^{+}_iq(q-q^{-1})^2.
\end{eqnarray*}
The inverse of this isomorphism is given by:
\begin{eqnarray*}
\label{eq:qliso1inv}
K^{\pm 1}_i &\mapsto & x^{\pm 1}_i,\\
\label{eq:qliso2inv}
e^-_i &\mapsto & y_i-x^{-1}_i, \\
\label{eq:qliso3inv}
e^+_i &\mapsto & (1-x_i z_i)q^{-1}(q-q^{-1})^{-2}.
\end{eqnarray*}
\end{proposition}

\begin{theorem}
{\rm \cite[Proposition 7.4]{qtet}}
\label{prop:uqhom}
For $i \in \Z_4$ there exists a $\C$-algebra homomorphism
from 
$U_q(\widehat{ \mathfrak{sl}}_2)$
to $\boxtimes_q$
that sends
\begin{eqnarray*}
&&
x_1 \mapsto x_{i,i+2},\quad
x^{-1}_1 \mapsto x_{i+2,i},\quad
y_1\mapsto x_{i+2,i+3},\quad
z_1 \mapsto x_{i+3,i},
\\
&&x_0 \mapsto x_{i+2,i}, \quad 
x^{-1}_0 \mapsto x_{i,i+2}, \quad
y_0 \mapsto x_{i,i+1},\quad
z_0 \mapsto x_{i+1,i+2}.
\end{eqnarray*}
\end{theorem}
\noindent {\it Proof:} 
Compare the defining relations for
$U_q(\widehat{ \mathfrak{sl}}_2)$
 given in
Proposition 
\ref{thm:qa2} with the relations
in Definition
\ref{def:qtet}.
\hfill $\Box $

\section{Distance-regular graphs; preliminaries}
We now turn our attention to distance-regular graphs.
After a brief review of the basic definitions
we recall the $Q$-polynomial property and the subconstituent
algebra.
For more information we refer the reader to 
\cite{bannai,bcn,godsil,terwSub1}.

\medskip
\noindent
Let $X$ denote a nonempty  finite  set.
Let $\hbox{Mat}_X(\C)$ 
denote the $\C$-algebra
consisting of all matrices whose rows and columns are indexed by $X$
and whose entries are in $\C  $. Let
$V=\C^X$ denote the vector space over $\C$
consisting of column vectors whose 
coordinates are indexed by $X$ and whose entries are
in $\C$.
We observe
$\hbox{Mat}_X(\C)$ 
acts on $V$ by left multiplication.
We call $V$ the {\it standard module}.
We endow $V$ with the Hermitean inner product $\langle \, , \, \rangle$ 
that satisfies
$\langle u,v \rangle = u^t\overline{v}$ for 
$u,v \in V$,
where $t$ denotes transpose and $\overline{\phantom{v}}$
denotes complex conjugation.
For all $y \in X,$ let $\hat{y}$ denote the element
of $V$ with a 1 in the $y$ coordinate and 0 in all other coordinates.
We observe $\{\hat{y}\;|\;y \in X\}$ is an orthonormal basis for $V.$

\medskip
\noindent
Let $\Gamma = (X,R)$ denote a finite, undirected, connected graph,
without loops or multiple edges, with vertex set $X$ and 
edge set
$R$.   
Let $\partial $ denote the
path-length distance function for $\Gamma $,  and set
$D := \mbox{max}\{\partial(x,y) \;|\; x,y \in X\}$.  
We call $D$  the {\it diameter} of $\Gamma $.
For an integer $k\geq 0$ we say that $\Gamma$ is {\it regular with
valency $k$} whenever each vertex of $\Gamma$ is adjacent to
exactly $k$ distinct vertices of $\Gamma$.
 We say that $\Gamma$ is {\it distance-regular}
whenever for all integers $h,i,j\;(0 \le h,i,j \le D)$ 
and for all
vertices $x,y \in X$ with $\partial(x,y)=h,$ the number
\begin{eqnarray*}
p_{ij}^h = |\{z \in X \; |\; \partial(x,z)=i, \partial(z,y)=j \}|
\end{eqnarray*}
is independent of $x$ and $y.$ The $p_{ij}^h$ are called
the {\it intersection numbers} of $\Gamma.$ 
We abbreviate $c_i=p^i_{1,i-1}$ $(1 \leq i \leq D)$,
$b_i=p^i_{1,i+1}$ $(0 \leq i \leq D-1)$,
$a_i=p^i_{1i}$ $(0 \leq i \leq D)$.

\medskip
\noindent
For the rest of this paper we assume  $\Gamma$  
is  distance-regular; to avoid trivialities we always
assume  $D\geq 3$.
Note that $\Gamma$ is regular with valency $k=b_0$. Moreover
 $k=c_i+a_i+b_i$ for $0 \leq i \leq D$, where $c_0=0$ and
$b_D=0$. 

\medskip
\noindent 
We mention a fact for later use.
By the triangle inequality, for $0 \leq h,i,j\leq D$ we have
$p^h_{ij}= 0$
(resp. 
$p^h_{ij}\not= 0$) whenever one of $h,i,j$ is greater than
(resp. equal to) the sum of the other two.

\medskip
\noindent 
We recall the Bose-Mesner algebra of $\Gamma.$ 
For 
$0 \le i \le D$ let $A_i$ denote the matrix in $\hbox{Mat}_X(\C)$ with
$(x,y)$-entry
$$
{(A_i)_{xy} = \cases{1, & if $\partial(x,y)=i$\cr
0, & if $\partial(x,y) \ne i$\cr}} \qquad (x,y \in X).
$$
We call $A_i$ the $i$th {\it distance matrix} of $\Gamma.$
We abbreviate $A=A_1$ and call this  the {\it adjacency
matrix} of $\Gamma.$ We observe
(i) $A_0 = I$;
 (ii)
$\sum_{i=0}^D A_i = J$;
(iii)
$\overline{A_i} = A_i \;(0 \le i \le D)$;
(iv) $A_i^t = A_i  \;(0 \le i \le D)$;
(v) $A_iA_j = \sum_{h=0}^D p_{ij}^h A_h \;( 0 \le i,j \le D)
$,
where $I$ (resp. $J$) denotes the identity matrix 
(resp. all 1's matrix) in 
 $\hbox{Mat}_X(\C)$.
 Using these facts  we find
 $\lbrace A_i\rbrace_{i=0}^D$
is a basis for a commutative subalgebra $M$ of 
$\mbox{Mat}_X(\C)$, called the 
{\it Bose-Mesner algebra} of $\Gamma$.
It turns out that $A$ generates $M$ \cite[p.~190]{bannai}.
By \cite[p.~45]{bcn}, $M$ has a second basis 
$\lbrace E_i\rbrace_{i=0}^D$ such that
(i) $E_0 = |X|^{-1}J$;
(ii) $\sum_{i=0}^D E_i = I$;
(iii) $\overline{E_i} = E_i \;(0 \le i \le D)$;
(iv) $E_i^t =E_i  \;(0 \le i \le D)$;
(v) $E_iE_j =\delta_{ij}E_i  \;(0 \le i,j \le D)$.
We call $\lbrace E_i\rbrace_{i=0}^D$  the {\it primitive idempotents}
of $\Gamma$.  

\medskip
\noindent
We  recall the eigenvalues
of  $\Gamma $.
Since $\lbrace E_i\rbrace_{i=0}^D$ form a basis for  
$M$ there exist complex scalars 
$\lbrace\theta_i\rbrace_{i=0}^D$
such that
$A = \sum_{i=0}^D \theta_iE_i$.
Observe
$AE_i = E_iA =  \theta_iE_i$ for $0 \leq i \leq D$.
By \cite[p.~197]{bannai} the 
scalars $\lbrace \theta_i\rbrace_{i=0}^D$ are
in $\R.$ Observe
$\lbrace\theta_i\rbrace_{i=0}^D$ are mutually distinct 
since $A$ generates $M$. We call $\theta_i$  the {\it eigenvalue}
of $\Gamma$ associated with $E_i$ $(0 \leq i \leq D)$.
Observe 
\begin{eqnarray*}
V = E_0V+E_1V+ \cdots +E_DV \qquad \qquad {\rm (orthogonal\ direct\ sum}).
\end{eqnarray*}
For $0 \le i \le D$ the space $E_iV$ is the  eigenspace of $A$ associated 
with $\theta_i$.

\medskip
\noindent 
We now recall the Krein parameters.
Let $\circ $ denote the entrywise product in
$\mbox{Mat}_X(\C)$.
Observe
$A_i\circ A_j= \delta_{ij}A_i$ for $0 \leq i,j\leq D$,
so
$M$ is closed under
$\circ$. Thus there exist complex scalars
$q^h_{ij}$  $(0 \leq h,i,j\leq D)$ such
that
$$
E_i\circ E_j = |X|^{-1}\sum_{h=0}^D q^h_{ij}E_h
\qquad (0 \leq i,j\leq D).
$$
By \cite[p.~170]{Biggs}, 
$q^h_{ij}$ is real and nonnegative  for $0 \leq h,i,j\leq D$.
The $q^h_{ij}$ are called the {\it Krein parameters} of $\Gamma$.
The graph $\Gamma$ is said to be {\it $Q$-polynomial}
(with respect to the given ordering
$\lbrace E_i\rbrace_{i=0}^D$
of the primitive idempotents)
whenever for $0 \leq h,i,j\leq D$, 
$q^h_{ij}= 0$
(resp. 
$q^h_{ij}\not= 0$) whenever one of $h,i,j$ is greater than
(resp. equal to) the sum of the other two
\cite[p.~235]{bcn}. See
\cite{
wdw,
caugh1,
caugh2,
curtin3,
curtin4,
dickie1,
dickie2,
aap1}
 for background information on the $Q$-polynomial property.
From now on we assume $\Gamma$ is $Q$-polynomial
with respect to $\lbrace E_i\rbrace_{i=0}^D$. We call the sequence
$\lbrace \theta_i\rbrace_{i=0}^D$ the {\it eigenvalue sequence}
for this $Q$-polynomial structure.

\medskip
\noindent
We  recall the dual Bose-Mesner algebra of $\Gamma.$
For the rest of this paper we fix
a vertex $x \in X.$ We view $x$ as a ``base vertex.''
For 
$ 0 \le i \le D$ let $E_i^*=E_i^*(x)$ denote the diagonal
matrix in $\hbox{Mat}_X(\C)$ with $(y,y)$-entry
\begin{equation}\label{DEFDEI}
{(E_i^*)_{yy} = \cases{1, & if $\partial(x,y)=i$\cr
0, & if $\partial(x,y) \ne i$\cr}} \qquad (y \in X).
\end{equation}
We call $E_i^*$ the  $i$th {\it dual idempotent} of $\Gamma$
 with respect to $x$ \cite[p.~378]{terwSub1}.
We observe
(i) $\sum_{i=0}^D E_i^*=I$;
(ii) $\overline{E_i^*} = E_i^*$ $(0 \le i \le D)$;
(iii) $E_i^{*t} = E_i^*$ $(0 \le i \le D)$;
(iv) $E_i^*E_j^* = \delta_{ij}E_i^* $ $(0 \le i,j \le D)$.
By these facts 
$\lbrace E_i^*\rbrace_{i=0}^D$ form a 
basis for a commutative subalgebra
$M^*=M^*(x)$ of 
$\hbox{Mat}_X(\C).$ 
We call 
$M^*$ the {\it dual Bose-Mesner algebra} of
$\Gamma$ with respect to $x$ \cite[p.~378]{terwSub1}.
For $0 \leq i \leq D$ let $A^*_i = A^*_i(x)$ denote the diagonal
matrix in 
 $\hbox{Mat}_X(\C)$
with $(y,y)$-entry
$(A^*_i)_{yy}=\vert X \vert (E_i)_{xy}$ for $y \in X$.
Then $\lbrace A^*_i\rbrace_{i=0}^D$ is a basis for $M^*$ 
\cite[p.~379]{terwSub1}.
Moreover
(i) $A^*_0 = I$;
(ii)
$\overline{A^*_i} = A^*_i \;(0 \le i \le D)$;
(iii) $A^{*t}_i = A^*_i  \;(0 \le i \le D)$;
(iv) $A^*_iA^*_j = \sum_{h=0}^D q_{ij}^h A^*_h \;( 0 \le i,j \le D)
$
\cite[p.~379]{terwSub1}.
We call 
 $\lbrace A^*_i\rbrace_{i=0}^D$
the {\it dual distance matrices} of $\Gamma$ with respect to $x$.
We abbreviate 
$A^*=A^*_1$ 
and call this the {\it dual adjacency matrix} of $\Gamma$ with
respect to $x$.
The matrix $A^*$ generates $M^*$ \cite[Lemma 3.11]{terwSub1}.

\medskip
\noindent We recall the dual eigenvalues of $\Gamma$.
Since $\lbrace E^*_i\rbrace_{i=0}^D$ form a basis for  
$M^*$ there exist complex scalars $\lbrace \theta^*_i\rbrace_{i=0}^D$
such that
$A^* = \sum_{i=0}^D \theta^*_iE^*_i$.
Observe
$A^*E^*_i = E^*_iA^* =  \theta^*_iE^*_i$ for $0 \leq i \leq D$.
By \cite[Lemma 3.11]{terwSub1} the 
scalars $\lbrace \theta^*_i\rbrace_{i=0}^D$ are in $\R$. 
The scalars $\lbrace \theta^*_i\rbrace_{i=0}^D$ are mutually
distinct 
since $A^*$ generates $M^*$. We call $\theta^*_i$ the {\it dual eigenvalue}
of $\Gamma$ associated with $E^*_i$ $(0 \leq i\leq D)$.
We call the sequence $\lbrace \theta^*_i\rbrace_{i=0}^D$ the
{\it dual eigenvalue sequence} for the given $Q$-polynomial structure.

\medskip
\noindent 
We recall the subconstituents of $\Gamma $.
From
(\ref{DEFDEI}) we find
\begin{equation}\label{DEIV}
E_i^*V = \mbox{span}\{\hat{y} \;|\; y \in X, \quad \partial(x,y)=i\}
\qquad (0 \le i \le D).
\end{equation}
By 
(\ref{DEIV})  and since
 $\{\hat{y}\;|\;y \in X\}$ is an orthonormal basis for $V$
 we find
\begin{eqnarray*}
\label{vsub}
V = E_0^*V+E_1^*V+ \cdots +E_D^*V \qquad \qquad 
{\rm (orthogonal\ direct\ sum}).
\end{eqnarray*}
For $0 \leq i \leq D$ the space $E^*_iV$ is the eigenspace
of $A^*$ associated with $\theta^*_i$.
We call $E_i^*V$ the $i$th {\it subconstituent} of $\Gamma$
with respect to $x$.

\medskip
\noindent
We recall the subconstituent algebra of $\Gamma $.
Let $T=T(x)$ denote the subalgebra of $\hbox{Mat}_X(\C)$ generated by 
$M$ and $M^*$. 
We call $T$ the {\it subconstituent algebra} 
(or {\it Terwilliger algebra}) of $\Gamma$ 
 with respect to $x$ \cite[Definition 3.3]{terwSub1}.
Observe that $T$ has finite dimension. Moreover $T$ is 
semisimple since it
is closed under the conjugate transponse map
\cite[p.~157]{CR}. We note that $A,A^*$ together generate $T$.
By \cite[Lemma 3.2]{terwSub1}
the following are relations in $T$:
\begin{eqnarray}
E^*_hA_iE^*_j&=&0 \quad \mbox{iff} \quad p^h_{ij}=0,
\qquad \qquad (0 \leq h,i,j \leq D),
\label{eq:triple2}
\\
E_hA^*_iE_j&=&0 \quad \mbox{iff} \quad q^h_{ij}=0, \qquad \qquad
(0 \leq h,i,j \leq D).
\label{eq:triple1}
\end{eqnarray}
See
\cite{curtin1,
curtin2,
curtin6,
egge1,
go,
go2,
hobart,
tanabe,
terwSub1,
terwSub2,
terwSub3}
for more information on the subconstituent
algebra.

\medskip
\noindent We recall the $T$-modules.
By a {\it T-module}
we mean a subspace $W \subseteq V$ such that $BW \subseteq W$
for all $B \in T.$ 
\noindent
Let $W$ denote a $T$-module and let 
$W'$ denote a  
$T$-module contained in $W$.
Then the orthogonal complement of $W'$ in $W$ is a $T$-module 
\cite[p.~802]{go2}.
It follows that each $T$-module
is an orthogonal direct sum of irreducible $T$-modules.
In particular $V$ is an orthogonal direct sum of irreducible $T$-modules.

\medskip
\noindent 
Let $W$ denote an irreducible $T$-module.
Observe that $W$ is the direct sum of the nonzero spaces among
$E^*_0W,\ldots, E^*_DW$. Similarly
$W$ is the direct sum 
 of the nonzero spaces among
$E_0W,\ldots,$ $ E_DW$.
By the {\it endpoint} of $W$ we mean
$\mbox{min}\lbrace i |0\leq i \leq D, \; E^*_iW\not=0\rbrace $.
By the {\it diameter} of $W$ we mean
$ |\lbrace i | 0 \leq i \leq D,\; E^*_iW\not=0 \rbrace |-1 $.
By the {\it dual endpoint} of $W$ we mean
$\mbox{min}\lbrace i |0\leq i \leq D, \; E_iW\not=0\rbrace $.
By
the {\it dual diameter} of $W$ we mean
$ |\lbrace i | 0 \leq i \leq D,\; E_iW\not=0 \rbrace |-1 $.
It turns out that the
diameter of $W$ is  equal to the dual diameter of
$W$
\cite[Corollary 3.3]{aap1}.
By \cite[Lemma 3.4]{terwSub1}
$\mbox{dim} \,E^*_iW \leq 1$ for $0 \leq i \leq D$ if and only if
$\mbox{dim} \,E_iW \leq 1$ for $0 \leq i \leq D$; in this case
$W$ is called {\it thin}. 

\medskip
\noindent
We finish this section with a few comments.

\begin{lemma}
{\rm \cite[Lemma 3.4, Lemma 3.9, Lemma 3.12]{terwSub1}}
\label{lem:basic}
Let $W$ denote an irreducible $T$-module with endpoint $\rho$,
dual endpoint $\tau$, and diameter $d$.
Then $\rho,\tau,d$ are nonnegative integers such that $\rho+d\leq D$ and
$\tau+d\leq D$. Moreover the following (i)--(iv) hold.
\begin{enumerate}
\item 
$E^*_iW \not=0$ if and only if $\rho \leq i \leq \rho+d$, 
$ \quad (0 \leq i \leq D)$.
\item
$W = \sum_{h=0}^{d} E^*_{\rho+h}W \qquad (\mbox{orthogonal direct sum}). $
\item 
$E_iW \not=0$ if and only if $\tau \leq i \leq \tau+d$,
$ \quad (0 \leq i \leq D)$.
\item
$W = \sum_{h=0}^{d} E_{\tau+h}W \qquad (\mbox{orthogonal direct sum}). $
\end{enumerate}
\end{lemma}

\begin{lemma}
\label{lem:ss}
{\rm 
\cite[Lemma~12.1]{qtetanddrg}}
For $Y \in 
 \hbox{\rm Mat}_X(\C)$ the following are equivalent:
\begin{enumerate}
\item $Y \in T$;
\item $YW\subseteq W$ for all irreducible $T$-modules $W$.
\end{enumerate}
\end{lemma}

\section{The split decomposition}

\noindent We are going to make use of a certain decomposition
of $V$ called the {\it split decomposition}.
The split decomposition was defined in \cite{ds} and discussed further in
\cite{qtetanddrg,jkim}. 
In this section we recall some results
on this topic.

\begin{definition}
\label{def:updown}
\rm
\cite[Definition~5.1]{ds}
For $-1\leq i,j\leq D$
we define
\begin{eqnarray*}
V_{i,j}^{\downarrow \downarrow} &=& 
(E^*_0V+\cdots+E^*_iV)\cap (E_0V+\cdots+E_jV),
\\
%V_{i,j}^{\uparrow \downarrow} &=& 
%(E^*_DV+\cdots+E^*_{D-i}V)\cap (E_0V+\cdots+E_jV).
%\\
V_{i,j}^{\downarrow \uparrow} &=& 
(E^*_0V+\cdots+E^*_iV)\cap (E_DV+\cdots+E_{D-j}V).
%\\
%V_{i,j}^{\uparrow \uparrow} &=& 
%(E^*_DV+\cdots+E^*_{D-i}V)\cap (E_DV+\cdots+E_{D-j}V).
\end{eqnarray*}
In the above two equations we interpret the
right-hand side to be 0 if $i=-1$ or $j=-1$. 
\end{definition}

\begin{definition}
\label{def:vtilde}
\rm
\cite[Definition~5.5]{ds}
With reference to Definition
\ref{def:updown}, for $(\mu,\nu)=
(\downarrow,\downarrow)$ or
$(\mu,\nu)=
(\downarrow,\uparrow)$ 
we have
$
V^{\mu \nu}_{i-1,j} \subseteq V^{\mu \nu}_{i,j}$
and
$
V^{\mu \nu}_{i,j-1}  \subseteq V^{\mu \nu}_{i,j}
$. Therefore
\begin{eqnarray*}
V^{\mu \nu}_{i-1,j}+
V^{\mu \nu}_{i,j-1} \subseteq V^{\mu \nu}_{i,j}.
\end{eqnarray*}
Referring to the above inclusion, we define ${\tilde V}^{\mu \nu}_{i,j}$
 to be the orthogonal complement of the left-hand side in the
 right-hand side; that
is
\begin{eqnarray*}
 {\tilde V}^{\mu \nu}_{i,j}=(
V^{\mu \nu}_{i-1,j}+
V^{\mu \nu}_{i,j-1})^\perp \cap V^{\mu \nu}_{i,j}.
\end{eqnarray*}
\end{definition}

\noindent The following is a mild generalization of
\cite[Corollary~5.8]{ds}.

\begin{lemma}
\label{lem:splitdec}
With reference to Definition
\ref{def:vtilde} the following holds 
for
$(\mu,\nu)=
(\downarrow,\downarrow)$ and 
$(\mu,\nu)=
(\downarrow,\uparrow)$:
\begin{eqnarray}
V = \sum_{i=0}^D \sum_{j=0}^D {\tilde V}^{\mu \nu}_{i,j}
\qquad \qquad (\mbox{\rm direct sum}).
\label{eq:spltdec}
\end{eqnarray}
\end{lemma}
\noindent {\it Proof:}
For $(\mu,\nu)=(\downarrow,\downarrow)$
 this is just \cite[Corollary~5.8]{ds}.
For $(\mu,\nu)=(\downarrow,\uparrow)$,
 in the proof
 of \cite[Corollary~5.8]{ds} replace the sequence
 $\lbrace E_i\rbrace_{i=0}^D$
 by $\lbrace E_{D-i}\rbrace_{i=0}^D$.
  \hfill $\Box $ \\

\begin{note} \rm
Following 
\cite[Definition~6.4]{jkim} 
we call the sum (\ref{eq:spltdec}) the 
 {\it 
 $(\mu,\nu)$-split
 decompostion of $V$}.
\end{note}

\noindent We now recall how the split decompositions are related
to the irreducible $T$-modules.
we start with a definition.

\begin{definition}
\label{displcement}
\rm 
\cite[Definition 4.1]{ds}
Let $W$ denote an irreducible $T$-module with endpoint
$\rho$, dual endpoint $\tau$, and diameter $d$. By the
{\it displacement of $W$ of the first kind} we mean
the scalar $\rho+\tau+d-D$. 
By the {\it displacement of $W$ of the second kind} we mean
the scalar $\rho-\tau$. 
By the inequalities in
Lemma \ref{lem:basic},
 each kind of displacement is at least $-D$ and at most
$D$.
\end{definition}

\begin{lemma}
\label{lem:twoview1}
{\rm \cite[Theorem 6.2]{ds}}
For $-D \leq \delta \leq D$ the following coincide:
\begin{enumerate}
\item The subspace of $V$ spanned by the irreducible $T$-modules
for which $\delta$ is the displacement of the first kind;
\item $\sum {\tilde V}^{\downarrow \downarrow}_{ij}$, where the
sum is over all ordered pairs $i,j$ $(0 \leq i,j\leq D)$
such that $i+j=\delta+D$.
\end{enumerate}
\end{lemma}

\begin{lemma}
\label{lem:twoview2}
For $-D \leq \delta \leq D$ the following coincide:
\begin{enumerate}
\item The subspace of $V$ spanned by the irreducible $T$-modules
for which $\delta$ is the displacement of the second kind;
\item $\sum {\tilde V}^{\downarrow \uparrow}_{ij}$, where the
sum is over all ordered pairs $i,j$ $(0 \leq i,j\leq D)$
such that $i+j=\delta+D$.
\end{enumerate}
\end{lemma}
\noindent {\it Proof:} 
In the proof of \cite[Theorem 6.2]{ds}, replace the
sequence $\lbrace E_i\rbrace_{i=0}^D$ by the
sequence $\lbrace E_{D-i}\rbrace_{i=0}^D$.
\hfill $\Box $ \\

\section{A homomorphism $\vartheta:\boxtimes_q \to T$}

\noindent We now impose an assumption on $\Gamma$.

\begin{assumption}
\label{def:sdcp}
\rm
We fix complex scalars $q,\eta,\eta^*,u,u^*,v,v^*$ with
$q,u,u^*,v,v^*$ nonzero 
and $q^{2i}\not=1$ for $1 \leq i \leq D$. We assume
that $\Gamma$ has a $Q$-polynomial structure with
eigenvalue sequence 
\begin{eqnarray*}
\theta_i &=& 
 \eta + uq^{2i-D}+vq^{D-2i} \qquad (0 \leq i \leq D)
\end{eqnarray*}
and dual eigenvalue sequence 
\begin{eqnarray*}
\theta^*_i &=& 
 \eta^* + u^*q^{2i-D}+v^*q^{D-2i} \qquad (0 \leq i \leq D).
\end{eqnarray*}
Moreover we assume that each irreducible $T$-module is thin.
\end{assumption}

\begin{remark} \rm
In the notation of Bannai and Ito \cite[p.~263]{bannai}
the $Q$-polynomial structure from 
Assumption 
\ref{def:sdcp}
is type I with $s\not=0, s^*\not=0$. We caution
the reader that the scalar denoted $q$ in
\cite[p.~263]{bannai} is the same as our scalar $q^2$.
\end{remark}

\begin{example}
\label{ex:nicegraph}
\rm The ordinary cycles are the only  
known distance-regular graphs that satisfy
 Assumption
\ref{def:sdcp}
\cite{bcn}.

\end{example}

\noindent
Under Assumption \ref{def:sdcp}
we will display a $\C$-algebra homomorphism
$\vartheta:\boxtimes_q \to T$.
To describe this
homomorphism we define two matrices in 
 $\hbox{Mat}_X(\C)$, called
$\Phi$ and $\Psi$.

\begin{definition}
\label{def:listdef}
\rm
With reference to 
Lemma \ref{lem:splitdec} and
Assumption \ref{def:sdcp},
let $\Phi$ (resp. $\Psi$) denote the unique 
 matrix
in
 $\hbox{Mat}_X(\C)$ 
that acts on 
${\tilde V}^{\downarrow \downarrow}_{ij}$ 
(resp. 
${\tilde V}^{\downarrow \uparrow}_{ij}$) 
as
$q^{i+j-D}I$ for $0 \leq i,j\leq D$.
Observe that each of $\Phi, \Psi$ is invertible.
\end{definition}

\begin{lemma}
\label{lem:phid}
Under Assumption \ref{def:sdcp} let
$W$ denote an irreducible $T$-module with endpoint $\rho$,
dual endpoint $\tau$, and diameter $d$. Then $\Phi$ and $\Psi$
act on
$W$ as $q^{\rho+\tau+d-D}I$ and $q^{\rho-\tau}I$ respectively.
\end{lemma}
\noindent {\it Proof:}
Concerning $\Phi$, abbreviate $\delta=\rho+\tau+d-D$ and recall that this
is the displacement of $W$ of the first kind. We show 
that $\Phi$ acts on $W$ as $q^{\delta}I$.
Let $V_{\delta}$ denote the common subspace from
parts (i), (ii) of Lemma
\ref{lem:twoview1}. By
 Lemma
\ref{lem:twoview1}(i) we have $W\subseteq V_\delta$.
 In Lemma
\ref{lem:twoview1}(ii) $V_\delta$ is expressed as a sum.
The matrix $\Phi$ acts on each term of this sum as
 $q^{\delta}I$ by
Definition
\ref{def:listdef},
so $\Phi$ acts  on $V_{\delta}$ as $q^{\delta}I$.
By these comments $\Phi$ acts on $W$ as $q^{\delta}I$
and this proves our assertion concerning $\Phi$. Our assertion
concerning $\Psi$ is similarly proved using the displacement of
the second kind and Lemma
\ref{lem:twoview2}.
\hfill $\Box $ \\

\begin{lemma}
\label{lem:cent}
Under Assumption \ref{def:sdcp} the matrices
$\Phi$ and $\Psi$ are central elements of $T$.
\end{lemma}
\noindent {\it Proof:}
The matrices $\Phi$ and $\Psi$ are contained in
$T$ by
Lemma \ref{lem:ss} and
 Lemma 
\ref{lem:phid}.
These 
matrices
are central in $T$ since by
Lemma 
\ref{lem:phid}
they act as a scalar multiple
of the identity on every irreducible $T$-module.
\hfill $\Box $ \\

\noindent The following is our main result.

\begin{theorem}
\label{thm:mres}
Under Assumption \ref{def:sdcp} there exists 
a $\C$-algebra homomorphism $\vartheta :\boxtimes_q \to T$
such that both
\begin{eqnarray}
A &=& \eta I + u \Phi \Psi^{-1}\vartheta(x_{01})+
v\Psi \Phi^{-1}\vartheta(x_{12}),
\label{eq:main1}
\\
A^* &=&\eta^* I + u^* \Phi \Psi \vartheta(x_{23})+
v^*\Psi^{-1} \Phi^{-1}\vartheta(x_{30}).
\label{eq:main2}
\end{eqnarray}
\end{theorem}

\noindent We will prove the above theorem after a few lemmas.

\begin{lemma}
\label{lem:localtet}
Under Assumption \ref{def:sdcp} let $W$
denote an irreducible $T$-module with endpoint $\rho$,
dual endpoint $\tau$, and diameter $d$. Then there exists
a $\boxtimes_q$-module structure on $W$ such that
the adjacency matrix $A$ 
acts as
$\eta I + u q^{2\tau+d-D}x_{01}+vq^{D-d-2\tau}x_{12}$
and the dual adjacency matrix $A^*$ acts as
$\eta^* I + u^* q^{2\rho+d-D}x_{23}+v^* q^{D-d-2\rho}x_{30}$.
This $\boxtimes_q$-module structure is irreducible.
\end{lemma}
\noindent {\it Proof:}
By \cite[Example~1.4]{itt} and since the $T$-module $W$ is thin the pair
$A,A^*$ acts on $W$ as a Leonard pair in the sense of
\cite[Definition~1.1]{qrac}. In the notation of \cite[Definition~5.1]{qrac}
this Leonard pair has an eigenvalue sequence 
$\lbrace \theta_{\tau+i}\rbrace_{i=0}^d$ and a 
dual eigenvalue sequence 
$\lbrace \theta^*_{\rho+i}\rbrace_{i=0}^d$ in view of
Lemma
\ref{lem:basic}.
To motivate what follows we note that
\begin{eqnarray*}
\label{eq:adjth}
\theta_{\tau+i} &=& \eta + u q^{2\tau+d-D} q^{2i-d} + 
 v q^{D-d-2\tau} q^{d-2i},
\\
\label{eq:adjths}
\theta^*_{\rho+i} &=& \eta^* + u^* q^{2\rho+d-D} q^{2i-d} + 
 v^* q^{D-d-2\rho} q^{d-2i}
\end{eqnarray*}
for $0 \leq i\leq d$. In both equations above the coefficients
of $q^{2i-d}$ and $q^{d-2i}$ are nonzero; therefore the action
of $A,A^*$ on $W$ is a Leonard pair of $q$-Racah type in the
sense of 
\cite[Example 5.3]{TLT:array}. Referring to this Leonard pair,
let $\lbrace \varphi_i\rbrace_{i=1}^d$
(resp. 
$\lbrace \phi_i\rbrace_{i=1}^d$) 
denote the first (resp. second) split sequence \cite[Section~7]{qrac}
associated with 
the eigenvalue sequence 
$\lbrace \theta_{\tau+i}\rbrace_{i=0}^d$ 
and  
the dual eigenvalue sequence 
$\lbrace \theta^*_{\rho+i}\rbrace_{i=0}^d$.
 By
\cite[Section~7]{qrac}
each of $\varphi_i,\phi_i$ is nonzero for $1 \leq i\leq d$.
By 
\cite[Example 5.3]{TLT:array} there exists a nonzero $r \in \C$ such
that 
\begin{eqnarray*}
&&\varphi_i = (q^i-q^{-i})(q^{d-i+1}-q^{i-d-1}) \\
&&  \qquad \qquad \times \quad (q^{d-i}-r^{-1}q^{i-1})
(uu^*rq^{2\tau+2\rho+d+i-2D}-vv^*q^{2D-2d-2\tau-2\rho+1-i}),
\\
&&\phi_i = (q^i-q^{-i})(q^{d-i+1}-q^{i-d-1}) \\
&&\qquad \qquad \times \quad (urq^{2\tau+d-D+1-i}-vq^{D-2d-2\tau+i})
(u^*q^{2\rho+d-D+i-1}-v^*r^{-1}q^{D-2\rho-i})
\end{eqnarray*}
for $1 \leq i \leq d$. 
Observe that
 $r$ is not among
$q^{d-1}, q^{d-3}, \ldots, q^{1-d}$ since
each of $\varphi_1, \varphi_2,\ldots,\varphi_d$ is nonzero.
By \cite[Section~7]{qrac} there exists a basis $\lbrace v_i\rbrace_{i=0}^d$
of $W$ such that
\begin{eqnarray*}
Av_i &=& \theta_{\tau+d-i}v_i+v_{i+1} \qquad (0 \leq i \leq d-1),
\quad Av_d=\theta_{\tau}v_d,
\\
A^*v_i &=& \theta^*_{\rho+i}v_i+\phi_i v_{i-1} \qquad (1 \leq i \leq d),
\quad A^*v_0=\theta^*_{\rho}v_0.
\end{eqnarray*}
For convenience we adjust this basis slightly.
For $1 \leq i \leq d$ define 
\begin{eqnarray*}
\gamma_i = 
(q^i-q^{-i})
(urq^{2\tau+d-D+1-i}-vq^{D-2d-2\tau+i}).
\end{eqnarray*}
In the above equation the right-hand side is nonzero
since it is a factor
of $\phi_i$, so
$\gamma_i \not=0$.
Define $u_i=(\gamma_1\gamma_2\cdots \gamma_i)^{-1}v_i$ 
for $0 \leq i \leq d$ 
and note that
 $\lbrace u_i\rbrace_{i=0}^d$ 
is a basis for $W$. 
By the construction
\begin{eqnarray*}
Au_i &=& \theta_{\tau+d-i}u_i+\gamma_{i+1} u_{i+1}
\qquad (0 \leq i \leq d-1),
\quad Au_d=\theta_{\tau}u_d,
\\
A^*u_i &=& \theta^*_{\rho+i}u_i+\phi_i \gamma^{-1}_i u_{i-1} 
\qquad (1 \leq i \leq d),
\quad A^*u_0=\theta^*_{\rho}u_0.
\end{eqnarray*}
We let each standard generator of $\boxtimes_q$
act linearly on $W$; to define this action we
specify 
what it does to
the basis $\lbrace u_i\rbrace_{i=0}^d$.
Here are the details:
\begin{eqnarray*}
&&x_{01}. u_i = q^{d-2i}u_i + (q^d-q^{d-2i-2})q^{1-d}ru_{i+1} \qquad
(0 \leq i \leq d-1), \qquad 
x_{01}. u_d = q^{-d}u_d,
\\
&&x_{12}. u_i = q^{2i-d}u_i + (q^{-d}-q^{2i+2-d})u_{i+1} \qquad
(0 \leq i \leq d-1), \qquad 
x_{12}. u_d = q^{d}u_d,
\\
&&x_{23}. u_i = q^{2i-d}u_i + (q^{d}-q^{2i-2-d})u_{i-1} \qquad
(1 \leq i \leq d), \qquad 
x_{23}. u_0 = q^{-d}u_0,
\\
&&x_{30}. u_i = q^{d-2i}u_i + (q^{-d}-q^{d-2i+2})q^{d-1}r^{-1}u_{i-1} \qquad
(1 \leq i \leq d), \qquad 
x_{30}. u_0 = q^{d}u_0,
\\
&&x_{13}.u_i = q^{2i-d}u_i \qquad (0 \leq i \leq d),
\\
&&x_{31}.u_i = q^{d-2i}u_i \qquad (0 \leq i \leq d),
\\
&&x_{02}.u_i =
(1-rq^{-d-1})\frac{(1-q^{2d-2i+2})(1-q^{2d-2i+4})\cdots (1-q^{2d})q^{d-2i}}
{(1-rq^{d-1-2i})(1-rq^{d+1-2i})\cdots
(1-rq^{d-1})} u_0
\\ 
&& \quad
 + \quad
(1-rq^{d+1})(1-rq^{-d-1})
\sum_{h=1}^i \frac{(1-q^{2d-2i+2})(1-q^{2d-2i+4})\cdots (1-q^{2d-2h})q^{d-2i}}
{(1-rq^{d-1-2i})(1-rq^{d+1-2i})\cdots
(1-rq^{d+1-2h})} u_h
\\
&& \quad
 + \quad
\frac{(q^{2i+2}-1)r}{q^{2i+1}(1-rq^{d-1-2i})}u_{i+1}
\qquad (0 \leq i \leq d-1),
\\
&&x_{02}.u_d =
\frac{(1-q^{2})(1-q^{4})\cdots (1-q^{2d})q^{-d}}
{(1-rq^{1-d})(1-rq^{3-d})\cdots
(1-rq^{d-1})} u_0
\\ 
&& \quad
 + \quad
(1-rq^{d+1})\sum_{h=1}^d \frac{(1-q^{2})(1-q^{4})\cdots
(1-q^{2d-2h})q^{-d}}
{(1-rq^{1-d})(1-rq^{3-d})\cdots
(1-rq^{d+1-2h})} u_h,
\\
&&x_{20}.u_0 = 
 (1-rq^{d+1})
\sum_{h=0}^{d-1}\frac{(1-q^{2})(1-q^{4})\cdots (1-q^{2h})r^hq^{h-dh-d}}
{(1-rq^{1-d})(1-rq^{3-d})\cdots (1-rq^{2h-d+1})}u_h
\\
&&\qquad \qquad + \quad 
\frac{(1-q^{2})(1-q^{4})\cdots (1-q^{2d})r^dq^{-d^2}}
{(1-rq^{1-d})(1-rq^{3-d})\cdots (1-rq^{d-1})} u_d,
\\
&&x_{20}.u_i = \frac{q^d-q^{2i-2-d}}{1-rq^{2i-d-1}}u_{i-1}
\\
&&+ \quad
(1-rq^{d+1})(1-rq^{-d-1})
\sum_{h=i}^{d-1}\frac{(1-q^{2i+2})(1-q^{2i+4})\cdots (1-q^{2h})r^{h-i}q^{(d+1)i-(d-1)h-d}} 
{(1-rq^{2i-d-1})(1-rq^{2i-d+1})\cdots (1-rq^{2h-d+1})}u_h
\\
&&+ \quad 
 (1-rq^{-d-1})
\frac{(1-q^{2i+2})(1-q^{2i+4})\cdots (1-q^{2d})r^{d-i}q^{di+i-d^2}}
{(1-rq^{2i-d-1})(1-rq^{2i-d+1})\cdots (1-rq^{d-1})}u_d
\qquad (1 \leq i \leq d).
\end{eqnarray*}
In the above formulae the denominators are nonzero since
$r$ is not among $q^{d-1}, q^{d-3}, \ldots, q^{1-d}$.
One checks (or see \cite{qtet24}) that the above actions
satisfy the defining relations for $\boxtimes_q$ from
Definition
\ref{def:qtet}, so these actions induce a $\boxtimes_q$-module structure
on $W$. Comparing the action of 
$A$ (resp. $A^*$) on $\lbrace u_i\rbrace_{i=0}^d$
with the actions of $x_{01}, x_{12}$ (resp. 
$x_{23}, x_{30}$) on 
$\lbrace u_i\rbrace_{i=0}^d$ we find that both
\begin{eqnarray*}
A&=&\eta I + u q^{2\tau+d-D}x_{01}+vq^{D-d-2\tau}x_{12},
\\
A^*&=&\eta^* I + u^* q^{2\rho+d-D}x_{23}+v^*q^{D-d-2\rho}x_{30}
\end{eqnarray*}
on $W$.
By these equations and since the $T$-module $W$ is irreducible
we find the $\boxtimes_q$-module $W$ is irreducible.
The result follows.
\hfill $\Box $ \\

\begin{lemma}
\label{lem:localtet2}
Under Assumption \ref{def:sdcp} let $W$
denote an irreducible $T$-module and consider the 
$\boxtimes_q$-action on $W$ from Lemma
\ref{lem:localtet}. Then the following equations hold on $W$:
\begin{eqnarray*}
A&=&\eta I + u \Phi \Psi^{-1}x_{01}+v\Psi \Phi^{-1}x_{12},
\\
A^*&=&\eta^* I + u^* \Phi \Psi x_{23}+v^* \Psi^{-1}\Phi^{-1}x_{30}.
\end{eqnarray*}
\end{lemma}
\noindent {\it Proof:} 
Combine Lemma
\ref{lem:phid}
and Lemma
\ref{lem:localtet}.
\hfill $\Box $ \\

\noindent It is now a simple matter to prove Theorem
\ref{thm:mres}.

\medskip
\noindent 
{\it Proof of Theorem \ref{thm:mres}}:
We start with a comment. Let $W$ and $W'$ denote
irreducible $T$-modules, and consider the 
$\boxtimes_q$-module structure on $W$ and $W'$
from Lemma
\ref{lem:localtet}. From the construction we may assume
that if the $T$-modules  
$W$ and $W'$ are isomorphic then the $\boxtimes_q$-modules 
$W$ and $W'$ are isomorphic. With our comment out of the
way we proceed to the main argument.
The standard module $V$ decomposes into a direct sum of
irreducible $T$-modules. Each irreducible $T$-module
in this decomposition supports a $\boxtimes_q$-module
structure from Lemma
\ref{lem:localtet}.
Combining these $\boxtimes_q$-modules we get a $\boxtimes_q$-module
structure on $V$.
This module structure induces
a $\C$-algebra homomorphism 
$\vartheta: \boxtimes_q \to 
\hbox{Mat}_X(\C)$.
The map $\vartheta$ satisfies
(\ref{eq:main1}), 
(\ref{eq:main2}) 
in view of  Lemma
\ref{lem:localtet2}.
To finish the proof it suffices to
show that $\vartheta(\boxtimes_q) \subseteq T$.
To this end we pick $\zeta \in \boxtimes_q$ 
and show
$\vartheta(\zeta)\in T$. 
Since $T$ is semisimple, and by our preliminary
comment, there exists $B \in T$
that acts on each irreducible
$T$-module in the above decomposition
as $\vartheta(\zeta)$. 
The $T$-modules in this decomposition
span
$V$ so $\vartheta(\zeta)$ coincides with $B$ on
$V$; therefore  $\vartheta(\zeta)=B$ and in particular
$\vartheta(\zeta) \in T$ as desired.
We have now shown that
$\vartheta(\boxtimes_q) \subseteq T$ and
the result follows.
\hfill $\Box $ \\

\begin{remark}
\label{rem:cp}
\rm 
In Theorem 
\ref{thm:mres} we obtained a $\C$-algebra homomorphism
 $\vartheta: \boxtimes_q\to T$.
In Theorem
\ref{prop:uqhom} we displayed four $\C$-algebra homomorphisms
from 
$U_q({\widehat{\mathfrak{sl}}_2})$ 
into $\boxtimes_q$. Composing these homomorphisms with
 $\vartheta $
we obtain four 
$\C$-algebra homomorphisms from 
$U_q({\widehat{\mathfrak{sl}}_2})$ into $T$.
\end{remark}

\noindent We conjecture that the conclusion of
Theorem 
\ref{thm:mres} still holds if we weaken Assumption
\ref{def:sdcp} by no longer requiring that each
irreducible $T$-module is thin.

\noindent Tatsuro Ito \hfil\break
\noindent Department of Computational Science \hfil\break
\noindent Faculty of Science \hfil\break
\noindent Kanazawa University \hfil\break
\noindent Kakuma-machi \hfil\break
\noindent Kanazawa 920-1192, Japan \hfil\break
\noindent email: {\tt tatsuro@kenroku.kanazawa-u.ac.jp}

\bigskip

\noindent Paul Terwilliger \hfil\break
\noindent Department of Mathematics \hfil\break
\noindent University of Wisconsin \hfil\break
\noindent 480 Lincoln Drive \hfil\break
\noindent Madison, WI 53706-1388 USA \hfil\break
\noindent email: {\tt terwilli@math.wisc.edu }\hfil\break


\begin{thebibliography}{10}

\bibitem{bannai} E. Bannai and T. Ito. {\it Algebraic Combinatorics I:
      Association Schemes.} Benjamin/Cummings, London, 1984.

%%%%%%%%%
%\bibitem{BT}
%G.~Benkart and P.~Terwilliger.
%\newblock The universal central extension of the three-point
%$\mathfrak{sl}_2$ loop algebra.
%\newblock{\em Proc. Amer. Math. Soc.}, accepted.
%{\tt arXiv:math.RA/0512422}.
%%%%%%%%%%%%%%

\bibitem{Biggs} N. Biggs {\it Algebraic Graph Theory. Second edition.}
      Cambridge University Press, Cambridge, 1993.

%%%%%%%%%
%\bibitem{Br}
%M.~Bremner.
%\newblock Generalized affine Kac-Moody Lie algebras over
%localizations of the polynomial ring in one variable.
%\newblock{\em Canad. Math. Bull.} {\bf 37} (1994), 21--28.
%%%%%%%%%%%

\bibitem{bcn} A. E. Brouwer, A. M. Cohen, and A. Neumaier. 
      {\it Distance-Regular Graphs.} Springer-Verlag, Berlin, 1989.


\bibitem{wdw} A. E. Brouwer, C. D. Godsil, J. H. Koolen,
W. J. Martin.
 Width and dual width of subsets in polynomial association
schemes.
\newblock{\em J. Combin. Theory Ser. A} {\bf 102}
(2003), 255-271.

\bibitem{caugh1} J. S. Caughman IV.
 Spectra of  bipartite $P$- and $Q$-polynomial association schemes.
\newblock{\em
Graphs Combin.} {\bf 14} (1998), 321--343.

\bibitem{caugh2} J. S. Caughman IV. 
 The Terwilliger algebras of bipartite $P$- and $Q$-polynomial
association
schemes.
\newblock{\em
Discrete Math.} {\bf 196} (1999), 65--95. 

\bibitem{charp}
V.~Chari and A.~Pressley.
\newblock {Q}uantum affine algebras.
\newblock {\em Comm. Math. Phys.}
{\bf 142} (1991), 261--283. 

\bibitem{curtin1} B. Curtin.
 Bipartite distance-regular graphs I.
\newblock{\em
Graphs Combin.} {\bf 15} (1999), 
143--158.

\bibitem{curtin2} B. Curtin. 
 Bipartite distance-regular graphs II.
\newblock{\em
Graphs Combin.} {\bf 15} (1999), 
377--391.

\bibitem{curtin3} B. Curtin.
 2-homogeneous bipartite distance-regular graphs.
\newblock{\em
Discrete Math.}
{\bf 187} (1998), 39--70.

\bibitem{curtin4} B. Curtin. 
 Distance-regular graphs which support a spin model are 
thin. 
\newblock{\em  16th British Combinatorial Conference (London, 1997).} Discrete Math.
{\bf 197/198} (1999), 205--216.

\bibitem{curtin6} B. Curtin and K. Nomura.
 Distance-regular graphs related to the quantum enveloping 
algebra of ${\rm sl}(2)$.
\newblock{\em
J. Algebraic Combin.}
{\bf 12} (2000), 25--36.

\bibitem{CR} C. Curtis and I. Reiner.
{\it {R}epresentation {T}heory of {F}inite {G}roups and
{A}ssociative {A}lgebras}.
Interscience, New York, 1962.
%%%%%%%%%


\bibitem{dickie1} G. Dickie. 
 Twice $Q$-polynomial distance-regular graphs 
 are thin.
\newblock{\em
European J. Combin.}
{\bf 16} (1995), 555--560. 

\bibitem{dickie2} G. Dickie and P. Terwilliger.
 A note on thin $P$-polynomial and dual-thin
$Q$-polynomial symmetric association  schemes.
\newblock{\em
J. Algebraic Combin.} {\bf 7} (1998), 5--15. 

%\bibitem{ding}
%J.~Ding and I.~B.~Frenkel.
%\newblock Isomorphism of two realizations of quantum affine
%algebra $U_q(gl(n))$.
%\newblock{\em Comm. Math. Phys.} {\bf 156} (1993), 277--300.
%%%%%%%%%%%%

\bibitem{egge1} E. Egge.
 A generalization of the Terwilliger algebra.
\newblock{\em
J. Algebra } {\bf 233} (2000), 213--252. 


%\bibitem{E}
%A.~Elduque.
%\newblock The $S_4$-action on the tetrahedron algebra.
%\newblock Preprint.  \hfil\break
%{\tt arXiv:math.RA/0604218}
%%%%%%%%%%%%%

\bibitem{go} J. T.  Go.
 The Terwilliger algebra of the hypercube.
\newblock{\em
European J. Combin.}
{\bf 23} (2002), 399--429.


\bibitem{go2} J. T. Go and P. Terwilliger.
 Tight distance-regular graphs and the subconstituent algebra.
\newblock{\em
 European J. Combin.} {\bf 23} (2002),
   793--816. 

\bibitem{godsil} C. D. Godsil. {\it Algebraic Combinatorics.}
Chapman and Hall, Inc., New York, 1993. 

%\bibitem{Ha}
%B.~Hartwig.
%\newblock The Tetrahedron algebra and its finite-dimensional
%irreducible modules.
%\newblock {\em Linear Algebra Appl.}, submitted.
%{\tt arXiv:math.RT/0606197}
%%%%%%%%%%%

\bibitem{HT}
B.~Hartwig and P.~Terwilliger.
\newblock The Tetrahedron algebra, the Onsager algebra,
and the $\mathfrak{sl}_2$ loop algebra.
\newblock{\em J. Algebra} {\bf 308} (2007), 840--863.
{\tt arXiv:math-ph/0511004}.

\bibitem{hobart} S. A. Hobart and T. Ito.  The structure of
nonthin irreducible $T$-modules: ladder bases and classical parameters.
\newblock{\em 
J. Algebraic Combin.} {\bf 7} (1998), 53--75.


\bibitem{itt} T. Ito, K. Tanabe,  P. Terwilliger.
 Some algebra related to $P$- and $Q$-polynomial association
schemes.
\newblock{\em
Codes and Association Schemes (Piscataway NJ, 1999), 167--192, DIMACS
Ser. Discrete Math. Theoret. Comput. Sci.}
{\bf 56}, Amer. Math. Soc., Providence RI 2001.
{\tt arXiv:math.CO/0406556}.


\bibitem{tdanduq}
T.~Ito and P.~Terwilliger.
\newblock {Tridiagonal pairs and the quantum affine 
algebra
$U_q({\widehat{sl}}_2)$.}
\newblock {\em Ramanujan J.}
{\bf 13} (2007), 39--62;
{\tt arXiv:math.QA/0310042}.

%\bibitem{NN}
%T.~Ito and P.~Terwilliger.
%\newblock Two non-nilpotent linear transformations that
%satisfy the cubic $q$-Serre relations.
%\newblock{\em J. Algebra Appl.}, submitted.
%{\tt arXiv:math.QA/0508398}.

%\bibitem{equit1}
%T.~Ito, P.~Terwilliger, and C.~W.~Weng.
%\newblock The quantum algebra 
%$U_q(\mathfrak{sl}_2)$ and its equitable presentation.
%\newblock {\em J. Algebra} {\bf 298} (2006), 284--301.
%{\tt arXiv:math.QA/0507477}. 

\bibitem{qtet}
T.~Ito and P.~Terwilliger.
\newblock The $q$-tetrahedron algebra and
its finite-dimensional irreducible modules.
\newblock {\em Comm. Algebra}; in press. {\tt arXiv:math.QA/0602199}.

\bibitem{qinv}
T.~Ito and P.~Terwilliger.
\newblock $q$-Inverting pairs of linear transformations
and the $q$-tetrahedron algebra.
\newblock {\em Linear Algebra Appl.}; in press.
{\tt arXiv:math.RT/0606237}. 

\bibitem{qtetanddrg}
T.~Ito and P.~Terwilliger.
\newblock Distance-regular graphs and the $q$-tetrahedron
algebra.
\newblock {\em European J. Combin.}; submitted. 
{\tt arXiv:math.CO/0608694}.

\bibitem{qtet24}
T.~Ito and P.~Terwilliger.
Evaluation modules for the $q$-tetrahedron algebra.
\newblock In preparation.



%%%%%%%%%need?
%\bibitem{kac}
%V.~G.~Kac.
%\newblock {\em Infinite dimensional Lie algebras}.
%Third Ed., 
%\newblock Cambridge U. Press, Cambridge, 1990.

%\bibitem{Kassel}
%C.~Kassel.
%\newblock {\em Quantum Groups},
%\newblock Springer-Verlag,
%New York,
%    1995.

\bibitem{jkim}
Joohyung~Kim.
A duality between pairs of split decompositions for
a $Q$-polynomial distance-regular graph.
\newblock {\em Discrete Math}; submitted.


\bibitem{lusztig}
G.~Lusztig.
\newblock  {\em Introduction to Quantum Groups},
\newblock  Birkhauser, Boston,
1990. 


\bibitem{aap1}
A. A. Pascasio. 
 On the multiplicities of the primitive idempotents of a $Q$-polynomial distance-regular graph.
\newblock{\em
European J.
   Combin.} {\bf 23} (2002),  1073--1078. 
    
%\bibitem{aap12} A. A. Pascasio. 
%Tight distance-regular
%graphs and the $Q$-polynomial property.
%\newblock{\em
%Graphs Combin. } {\bf 17} (2001), 149--169.

%\bibitem{PT}
%A.~A.~Pascasio and P.~Terwilliger.
%\newblock The Tetrahedron algebra and the Hamming graphs.
%\newblock In preparation.


\bibitem{tanabe} K. Tanabe. 
The irreducible modules of the Terwilliger
algebras of Doob schemes.
\newblock{\em
J. Algebraic Combin.} {\bf 6} (1997), 173--195.

%%%%%%%%%%%%%
%\bibitem{tarasov}
%V.~O.~Tarasov.
%\newblock Irreducible monodromy matrices for the $R$-matrix of the
%XXZ-model, and lattice local quantum Hamiltonians.
%\newblock{\em
%Theor. Math. Phys.} {\bf 63}
%(1985), 440--454.
%%%%%%%%%%%%%%


\bibitem{terwSub1} P. Terwilliger. 
The subconstituent algebra of
an association scheme I.
\newblock{\em
J. Algebraic Combin.}
{\bf 1} (1992), 363--388.  

\bibitem{terwSub2} P. Terwilliger. 
The subconstituent algebra of
an association scheme II.
\newblock{\em
J. Algebraic Combin.}
{\bf 2} (1993), 73--103.  

\bibitem{terwSub3} P. Terwilliger.  The subconstituent algebra of
an association scheme III. 
\newblock{ \em
J. Algebraic Combin.}
{\bf 2} (1993), 177--210.  


%\bibitem{LS99}
%P.~Terwilliger.
%\newblock Two linear transformations each tridiagonal with respect to an
% eigenbasis of the other.
%\newblock {\em Linear Algebra Appl.}  {\bf 330} (2001), 149--203.
		   
%\bibitem{LS24}
%P.~Terwilliger.
%\newblock  Leonard pairs from 24 points of view.
%\newblock {\em Rocky Mountain J. Math.} {\bf 32} (2002), 827--888.
			       
%\bibitem{lsint}
%P.~Terwilliger.
%\newblock Introduction to {L}eonard pairs.
%\newblock {OPSFA Rome 2001}.
%\newblock{\em J. Comput. Appl. Math.} {\bf 153}(2) (2003),
%463--475.
							     
\bibitem{qrac}
P. Terwilliger.
\newblock Leonard pairs and the $q$-Racah polynomials.
\newblock {\em Linear Algebra Appl.} {\bf 387} (2004), 235--276. 
{\tt arXiv:math.QA/0306301}.

\bibitem{TLT:array}
P.~Terwilliger.
\newblock Two linear transformations each tridiagonal with respect to an
eigenbasis of the other; comments on the parameter array.
\newblock{\em Des. Codes Cryptogr.} {\bf 34} (2005) 307--332.
{\tt arXiv:math.RA/0306291}.
												 
\bibitem{ds}
P.~Terwilliger.
\newblock The displacement and split decompositions
for a $Q$-polynomial distance-regular graph.
\newblock{\em Graphs Combin.} {\bf 21} (2005), 263--276.
{\tt arXiv:math.CO/0306142}.


\bibitem{equit2}
P.~Terwilliger.
\newblock 
The equitable presentation for the quantum group
$U_q({\mathfrak{g}})$ associated with a symmetrizable
Kac-Moody algebra $\mathfrak{g}$.
\newblock {\em J. Algebra} {\bf 298} (2006), 302--319.
{\tt arXiv:math.QA/0507478}.

%\bibitem{uniform}
%P.~Terwilliger.
%\newblock 
%The incidence algebra of a uniform poset. 
%\newblock{\em Coding theory and design theory, Part I},  193--212,
%\newblock IMA Vol. Math. Appl., {\bf 20}, Springer, New York, 1990.


 \end{thebibliography}
\end{document}